\documentclass[11pt, reqno]{amsart}
\usepackage{bm}
\usepackage{verbatim}

\usepackage{amsmath,amssymb,amscd,amsthm,amsxtra, esint, bbm}
\usepackage[dvips]{graphics,epsfig}
\usepackage{color}

\allowdisplaybreaks[2]

\definecolor{orange}{rgb}{1,0.5,0}
\definecolor{darkgreen}{rgb}{.2,0.4,0.3}

\newtheorem{theorem}{Theorem} [section]

\newtheorem{lemma}[theorem]{Lemma}
\newtheorem{proposition}[theorem]{Proposition}
\newtheorem{remark}[theorem]{Remark}

\newtheorem{corollary}[theorem]{Corollary}

\newcommand{\I}{\hspace{0.5mm}\textup{I}\hspace{0.5mm}}
\newcommand{\II}{\textup{I \hspace{-2.8mm} I} }

\newcommand{\noi}{\noindent}
\newcommand{\Z}{\mathbb{Z}}
\newcommand{\R}{\mathbb{R}}

\newcommand{\T}{\mathbb{T}}

\newcommand{\dl}{\delta}
\newcommand{\Dl}{\Delta}
\newcommand{\eps}{\varepsilon}

\newcommand{\G}{\Gamma}
\newcommand{\ld}{\lambda}

\newcommand{\s}{\sigma}
\newcommand{\ft}{\hat}
\newcommand{\wt}{\widetilde}
\newcommand{\cj}{\overline}

\newcommand{\dt}{\partial_t}

\newcommand{\betat}{{\tilde \beta}}

\newcommand{\too}{\longrightarrow}

\newcommand{\jb}[1]
{\langle #1 \rangle}

\newcommand{\proj}{\mathbf{P}}

\renewcommand{\b}{\beta}
\newcommand{\ind}{\mathbf 1}

\numberwithin{equation}{section}
\numberwithin{theorem}{section}

\title[Interpolation of Gibbs measures with White Noise for Hamiltonian PDE]
{Interpolation of Gibbs measures with White Noise for Hamiltonian PDE}

\author{Tadahiro Oh}
\address{\noindent Department of Mathematics, University of Toronto
\newline e-mail:  \rm \texttt{oh@math.toronto.edu}}
\author{Jeremy Quastel}
\address{\noindent Departments of Mathematics and Statistics, University of Toronto
\newline e-mail:  \rm \texttt{quastel@math.toronto.edu}}
\thanks{J.~Quastel was partially supported by the Natural Sciences and Engineering Research Council of Canada.}

\author{Benedek Valk\' o}
\address{\noindent Department of Mathematics, University of Wisconsin
\newline e-mail:  \rm \texttt{valko@math.wisc.edu }}
\thanks{B.~Valk\' o was  supported by the NSF Grant DMS-09-05820}

\begin{document}

\begin{abstract}  We consider the family of 
interpolation measures of Gibbs measures and white noise given by
$$dQ_{0,\b}^{(p)} =  Z_\b^{-1} \ind_{\{\int_{\T} u^2\le K\b^{-1/2}\}} e^{ -\int_{\T} u^2 +\b \int u^p } dP_{0,\b}$$
where $P_{0, \b} $ is the Wiener measure on the circle, with 
variance $\beta^{-1}$, conditioned to have mean zero.  
It is shown that
as $\beta\to 0$, $Q_0^\beta$ converges weakly to 
mean zero Gaussian white noise $Q_0$.  
As an application, we present  a straightforward proof that 
$Q_0$ is invariant for the Kortweg-de Vries equation (KdV).
This weak convergence also shows that the white noise 
is a weak limit of invariant measures for the modified KdV and the cubic nonlinear Schr\"odinger equations.
\end{abstract}



\date{\today}
\maketitle

\tableofcontents

\section{Introduction}

\subsection{An interpolation of measures}

Let $Q_0$ denote the mean zero Gaussian white noise on the circle $\T= \R/\Z$. 
i.e. $Q_0$ is the probability measure on real-valued distributions $u$ with $\int_{\mathbb T} u=0$ satisfying
\begin{equation}\label{white1}
\int e^{i \jb{f, u}} dQ_0(u) = e^{-\frac{1}{2} \|f\|_{L^2}^2}
\end{equation}

\noi
for any mean zero smooth real-valued function $f$ on $\mathbb T$,
where $\jb{\cdot, \cdot}$ denotes the pairing between the Schwartz space $\mathcal{S}(\T)$ and its dual $\mathcal{S}'(\T)$. 
It is known that 
$Q_0$ is supported on the Sobolev space $H^{s}_0(\T)$ for $s < -\frac{1}{2}$,
where $H^{s}_0(\T)$ consists of real-valued distributions 
$u =\sum_{n\neq 0} \hat u_n e^{2\pi i n x} \in \mathcal{S}'$  
with $\hat u_{-n}=\overline{\hat{u_n}}$ such that 
$ \|u\|_{H^{s}_0}^2 = \sum_{n\ne 0} |n|^{2s}|\hat u_n|^2 < \infty $.

Let $P_0$ denote the Wiener measure on $u\in C(\mathbb T)$ conditioned to have $\int_{\T} u= 0$.   
It can be derived from the Brownian Bridge $P$  as follows:  
For a given $x\in\R$,  condition a standard Brownian motion $u (t)$, $t\in [0,1]$, 
starting at $u(0)=x$ to have $u(1)=x$ and  $\int_{\T} u= 0$. Then distribute $u(0)$ according to a real Gaussian  with mean zero and variance $\pi^2/3$. The easiest way to check that this produces the appropriate measure is by the Fourier representation of $u$:
Let $\{g_n\}_{n\geq1}$ be a family of independent standard complex-valued Gaussian random variables, i.e. its real and imaginary parts are independent Gaussians with mean zero and variance $1/2$.   Also, for $n\geq 1$, let $g_{-n} = \cj{g_n}$.
Then 
 \begin{equation} \label{rep10}
u(x) = \sum_{n \ne 0} \frac{g_n}{n} e^{2\pi i n x}.\end{equation}
 
Similarly, let $P_{0, \b}$  be the Wiener measure with variance $\b^{-1}$ 
conditioned to have $\int_{\mathbb T} u= 0$. 
Formally, we can write $P_{0,\b}$ as
\begin{equation} \label{Wiener1}
 dP_{0,\beta} = Z_{0, \beta}^{-1} \exp\Big( -\frac{\b}{2}\int_{\T} u_x^2 \Big)\prod_{x\in \T} du(x).
\end{equation}
and under $P_{0, \b}$,  \begin{equation} \label{rep11}u(x) = \betat^{-1/2}\sum_{n \ne 0} \frac{g_n}{n} e^{2\pi i n x}, \quad \betat=4\pi^2 \beta.
\end{equation}
 
\noi 
For fixed $K>0$ and $p \in \mathbb{N}$,
let $P^{\varphi_1^p}_{0}$ denote the probability measure
on $u\in C(\T)$ with $\int_{\T} u=0$ given by
\begin{equation} \label{Gibbs1}
 dP^{\varphi_1^p}_{0} = Z_{p,K}^{-1}\ind_{\{\int_{\T} u^2\le K\}} e^{\int u^p } dP_0.
\end{equation}

\noi
The $L^2$-cutoff is necessary to make the normalization $Z_{p,K}$ well-defined and finite 
(for $ p \leq 6$ \cite{LRS, B2}.)
The notation $\varphi_1^p$ is borrowed from quantum field theory; 
the superscript $p$ denotes the order of the nonlinearity and the subscript the dimension. 
The measure $P^{\varphi_1^p}_{0}$ corresponds to the Gibbs measure for certain Hamiltonian PDEs. 
We will discuss this aspect in the next subsection. 
 
We can also define a family of probability measures depending on $\beta>0$,
\begin{equation} \label{Gibbs2}
dP^{\varphi_1^p}_{0,\b} = \hat Z_\b^{-1} \ind_{\{\int_{\T} u^2\le K\b^{-1/2}\}} e^{\b\int u^p } dP_{0, \b},
\end{equation}
 
\noi
where $\hat Z_\b = \hat Z (\b, p, K)$.
Finally, let  $Q^{p}_{0, \b}$, $\b>0$, be the following family of probability measures 
on $u\in C(\T)$ with $\int_{\T} u=0$, 
interpolating between $P^{\varphi_1^p}_{0,\beta}$ and $Q_0$;
\begin{equation} \label{Gibbs3}
dQ_{0, \b}^{(p)} =  Z_\b^{-1} \ind_{\{\int_{\T} u^2\le K\b^{-1/2}\}} e^{ -\frac{1}{2}\int_{\T} u^2 +\b \int u^p } dP_{0, \b}.
\end{equation}

\noi
In the following, we assume $p = 3$ or $4$.
It follows from \cite{LRS, B2}
that for each {\it fixed} $\b > 0$, 
$Q_{0, \b}^{(p)}$ is a well-defined probability measure on $H^s(\T)$, $s < \frac{1}{2}$, the regularity
being inherited from Brownian motion on $\T$.

The main result of this article is
\begin{theorem} \label{thm1}
Let $p = 3$ or $4$ and ${K>\frac12}$.
Then, as $\beta\to 0$, 
$Q_{0, \b}^{(p)}$ converges weakly to  $Q_0$ 
as probability measures on $H^{s}_0(\mathbb T)$, $s<-\frac{1}{2}$.
\end{theorem}

\begin{remark} \label{REM:complex} \rm
When $p = 4$, the analogue  to Theorem \ref{thm1} holds 
for the measures on complex-valued distributions $u$ (without the mean zero assumption),
\begin{equation} \label{Gibbs3c}
d{\bf Q}_{\b}^{(4)} =  \tilde Z_\b^{-1} \ind_{\{\int |u|^2\le K\b^{-1/2}\}} e^{ -\frac{1}{2}\int |u|^2 +\b \int |u|^4 } d{\bf P}_{\b},
\end{equation}

\noi
where ${\bf P}_{\b}$ is the complex  Wiener measure with variance $\b^{-1}$.
We present the proof of Theorem \ref{thm1} in details for the real-valued case
and indicate the modification for the complex-valued case.
\end{remark}

\noi
Formally, the theorem follows from the observation that
\begin{align} \label{conv1}
dQ_{0, \b}^{(p)} & = \bar{Z}_\b^{-1} \ind_{\{\int_{\T} u^2\le K\beta^{-1/2}\}} 
e^{ -\frac{1}{2} \int_{\T} u^2 +\b \int u^p -  \frac{\beta} 2 \int u_x^2 }\prod_{x\in \T} du(x) \\
& \stackrel{\beta\to 0}{\too} \quad\bar Z_0^{-1}e^{-\frac12\int_{\T} u^2}\prod_{x\in \T} du(x)= dQ_0.
 \notag
 \end{align}

\noi 
 So the result is intuitively clear.  Unfortunately, neither the normalizations $\bar Z_\b$ 
 nor the ``flat measure'' $\prod_{x\in \T} du(x)$ make sense, 
 so a proof is required.  
 It turns out to be a little tricky
 and it involves a careful analysis of random Fourier series.  

Consider the Gaussian measure $\mu_\beta$ given by 
\begin{align} \label{Gauss2}
d \mu_\beta 
= Z_\beta^{-1} e^{-\frac{1}{2} \int_{\T} u^2} d P_{0, \b}
=  \ft{Z}_\b^{-1} e^{-\frac{1}{2} \int_{\T} u^2 - \frac{\b}{2} \int u_x^2} \prod_{x \in \T} du(x), 
\end{align}

\noi
where $u$ is real-valued with $\int_\T u = 0$.\footnote{In the following, we use $Z_\b$
to denote various normalization constants.}
This is an interpolation of the Wiener measure $P_{0, \b}$ and the white noise $Q_0$ on $\T$.
{
If $u$ is distributed according to  $\mu_\beta$, then it can also  be represented as} 
\begin{equation} \label{rep1}
u(x) = \sum_{n \ne 0} \frac{g_n}{\sqrt{1  + \betat n^2}} e^{2\pi i n x}.\end{equation}

%

The main difficulty of the proof of Theorem \ref{thm1}
lies in establishing the exponential expectation estimate:
\begin{equation} \label{ExpEx1}
\mathbb{E}_{\mu_\beta} \big[\, \ind_{\{ \int_{\T} u^2 \leq K \b^{-\frac{1}{2}} \}} e^{r \b \int u^p}\big]
= \int \ind_{\{ \int_{\T} u^2 \leq K \b^{-\frac{1}{2}} \}} e^{r \b \int u^p} d \mu_\b
\leq C(r) <\infty
\end{equation}

\noi 
uniformly for small $\beta > 0$,
where $\mathbb{E}_{\mu_\beta}$ denotes an expectation with respect to $\mu_\b$.
Recall that for each $\beta > 0$, 
$u$ is almost surely in $H^s \setminus H^{\frac{1}{2}}$, $s <\frac{1}{2}$.
However, when $\beta = 0$, 
\eqref{Gauss2} reduces to the white noise $Q_0$
supported on $H^s \setminus H^{-\frac{1}{2}}$, $s<-\frac{1}{2}$. 
Hence, $\int u^p$, $p = 3, 4$,  diverges as $\beta \to 0$,
and thus we need to carefully analyze $\b \int u^p$ as $\beta \to 0$.
It turns out that the decay of $\b$ and the growth of $\int u^p$
is in perfect balance when $p = 4$, (see Remark \ref{REM:exp})
and the proof \eqref{ExpEx1} is much more delicate when $p = 4$.
We need some probabilistic tools such as
the hypercontractivity of the Ornstein-Uhlenbeck semigroup.
We present the proof in the remaining sections of the article.

\medskip
When $p = 4$, one can also consider the convergence of $\wt{Q}_{0, \b}^{(4)}$
whose density is given by
\[d \wt{Q}_{0, \b}^{(4)} = Z_\b^{-1} e^{- \beta \int_{\T} u^4} d \mu_\b.\]

\noi
In this case, thanks to the negative sign in front of $\beta \int_{\T} u^4$,
we have the exponential expectation estimate \eqref{ExpEx1} for free.
\begin{theorem} \label{thm0}
As $\beta\to 0$, 
$\wt{Q}_{0, \b}^{(4)}$ converges weakly to  $Q_0$ 
as probability measures on $H^{s}_0(\mathbb T)$, $s<-\frac{1}{2}$.
\end{theorem}

\noi
In proving Theorem \ref{thm0}, 
we follow the basic argument for Theorem \ref{thm1}.
However, since there is no need for an $L^2$-cutoff,
a slight care is required.
When $p = 3$, we still need an $L^2$-cutoff
in view of transformation $u \to -u$.

\medskip

Before we proceed to the proof of Theorem \ref{thm1}, 
let us discuss the motivation for studying this problem
and present an application to some Hamiltonian PDEs
in the remaining part of this section.

\subsection{Hamiltonian dynamics and Gibbs measures}
 Given a Hamiltonian flow on $\mathbb{R}^{2n}$:
\begin{equation} \label{HR2}
\begin{cases}
\dot{p}_i = \frac{\partial H}{\partial q_j}\\
\dot{q}_i = - \frac{\partial H}{\partial p_j} 
\end{cases}
\end{equation}

\noi
 with Hamiltonian $ H (p, q)= H(p_1, \dots, p_n, q_1, \dots, q_n)$,
Liouville's theorem states that the Lebesgue measure on $\mathbb{R}^{2n}$ is invariant under the flow.
Then, it follows from the conservation of the Hamiltonian $H$
that  the Gibbs measures $e^{-  H(p, q)} \prod_{j = 1}^{n} dp_j dq_j$ 
are invariant under the flow of \eqref{HR2}.

In the context of the nonlinear Schr\"odinger equations (NLS) on $\T$: 
\begin{equation} \label{NLS}
i u_t - u_{xx} \pm |u|^{p-2} u = 0, \qquad u(0) = 0,
\end{equation}

\noi
Lebowitz-Rose-Speer \cite{LRS} considered the Gibbs measure of the form 
\begin{equation} \label{Gibbs}
d \mu = \exp (- H(u)) \prod_{x\in \T} d u(x),
\end{equation}

\noi
where 
$H(u)$ is the Hamiltonian given by $H(u) = \frac{1}{2} \int |u_x|^2 \pm \frac{1}{p} \int |u|^p dx$.
It was shown that such Gibbs measure $\mu$ is  
a well-defined probability measure on $H^s \setminus H^{\frac{1}{2}}$, $s <\frac{1}{2}$.
(In the focusing case (with  $-$), 
the result only holds for $p < 6$ 
with the $L^2$-cutoff $\ind_{\{\int |u|^2 \leq K \}}$ for any $K>0$, 
and for $ p = 6$ with sufficiently small $K$.)
Using the Fourier analytic approach, Bourgain \cite{B2} continued the study and 
proved the invariance of the Gibbs measure $\mu$ under the flow of NLS.
In the same paper, he also established the invariance of the Gibbs measures
for the Korteweg-de Vries equation (KdV) on $\T$: 
\begin{equation}\label{KdV}
u_t +u_{xxx} -6 u u_x  =0, \qquad u(0) = u_0,
\end{equation}

\noi
and the modified KdV equation (mKdV) on $\T$:
\begin{equation}\label{mKdV}
u_t +u_{xxx} \mp  u^2 u_x  =0, \qquad u(0) = u_0.
\end{equation}

Invariant Gibbs measures $\mu$ for Hamiltonian PDEs can be regarded
as stationary measures for infinite dimensional dynamical systems,
and it follows from  Poincar\'e recurrence theorem 
that almost all the points of the phase space are stable according to Poisson, 
{i.e. if $\mathcal{S}_t$ denotes a flow map: $u_0 \mapsto u(t) = \mathcal{S}_{t} u_0$, 
then for almost all $u_0$, there exists a sequence $\{t_n\}$ tending to $\infty$
such that $ \mathcal{S}_{t_n} u_0 \to u_0$.}
We also know such dynamics is also multiply recurrent
in view of Furstenberg \cite{F}: 
{let $A$ be any measureable set with $\mu(A) > 0$.
Then, for any integer $k >1$, there exists $n \ne0$
such that 
$\mu( A \cap \mathcal{S}_n A
\cap \mathcal{S}_{2n} A \cap\cdots \cap\mathcal{S}_{(k-1)n} A )>0$.}
Note that this recurrence property holds only in the support of the Gibbs measure,
i.e. not for smooth functions.

Now note that if $F(p, q)$ is any function that is conserved under the flow of \eqref{HR2}, 
then  the measure $d \mu_F = e^{- F(p, q)} \prod_{j = 1}^{n} dp_j dq_j$ is invariant.
Recall that NLS, KdV, and mKdV are all Hamiltonian partial differential equations preserving the $L^2$-norm (see also  \cite{DLT} for another intriguing connection.)
 Hence, it is natural, at least at a heuristic level, 
to expect the invariance of the white noise for these equations.  
The difficulty here is the low regularity of the phase space.

\subsection{Invariance of white noise for KdV on $\T$.}

As an application of Theorem \ref{thm1}, 
we present a straightforward proof of the fact that
$Q_0$ is an invariant measure for KdV on $\T$.
Given a smooth initial condition $u_0:\mathbb T\to \mathbb R$, 
we have a solution $\mathcal{S}_t u_0 = u(t)$ for $-\infty<t < \infty$.  
In fact, KdV is well-posed for much rougher initial data;  
the nonlinear solution map $\mathcal S_t$ extends to 
a continuous group of nonlinear evolution operators
\begin{equation}\label{1point8}
\bar{\mathcal S}_t: H^{s}_0(\mathbb T) \to H^{s}_0(\T), \quad -\infty < t <\infty, \qquad s \geq -1.
\end{equation}

\noi
By the Fourier restriction method, 
Bourgain \cite{B1} proved $s \geq 0$, 
and Kenig-Ponce-Vega \cite{KPV} and {Colliander et al.}  \cite{CKSTT} pushed it down to $s \geq -\frac{1}{2}$.
Finally, Kappeler and Topalov \cite{KT} proved $s \geq -1$ via the inverse spectral method.
Since the white noise $Q_0$ is supported on $H_0^{s}(\T)$ for $s < -\frac{1}{2}$, 
this means that it makes sense to start KdV on the circle 
with white noise as initial data, for almost every realization.
 
In \cite{QV} and \cite{O1, O2}, we proved the following result:
\begin{theorem}\label{whiteinv}  
White noise $Q_0$ is invariant under KdV.  
i.e. for any $t\in \mathbb R$,
$\bar{\mathcal S}_t^{*} Q_0 = Q_0.$
\end{theorem}

\noi
Here, $\bar{\mathcal S}_t^{*}Q_0$ denotes the pushforward of the measure $Q_0$ by the map $\bar{\mathcal S}_t$.
The proof in \cite{QV} is indirect:  We show that
$Q_0$ is the image under the Miura transform of 
the Gibbs measure for the defocusing mKdV (with the $-$ sign in \eqref{mKdV}),
which was proven to be invariant by Bourgain \cite{B2}.
While the proof in \cite{O1, O2} is more direct,
it relies on heavy Fourier analysis.  Since the result is so simple to state, it is reasonable to ask for a straightforward proof (and 
such a proof has been requested of the authors.)

In the following, we give a {more straightforward proof} of Theorem \ref{whiteinv},
using Theorem \ref{thm1}, (\ref{1point8}), and the 
following.

\begin{proposition}   [Bourgain, \cite{B2}]  \label{thm2}
$P^{\varphi_1^3}_{0,\beta}$ defined in \eqref{Gibbs1}, $\b > 0$, are invariant 
for KdV.
\end{proposition}

\noi
Note that in \cite{B2} this is only explicitly proven for $\b=1$.  
But the same proof works for all $\b>0$.
If $\mu$ is an invariant measure of a Markov process $u(t)$ 
and $F$ is a conserved quantity; $F(u(t))=F(u(0))$, 
then, as long as it makes sense, $d\nu= Fd\mu$ is an invariant measure as well.
The quantity $F(u)= \int_{\T} u^2$ is a conserved quantity for KdV 
and $\exp(-\frac{1}{2}\int_{\T} u^2)\in L^1(P^{\varphi_1^3}_{0,\beta})$.
Hence it follows from Proposition \ref{thm2} that
\begin{corollary}\label{cor1}
$Q_{0, \b}^{(3)} $ defined in \eqref{Gibbs3}, $\beta>0$, are invariant for KdV.
\end{corollary}

\noi
To complete the proof of Theorem \ref{whiteinv}, 
we need to verify that $Q_0$, the limit
of invariant measures by Theorem \ref{thm1} and Corollary \ref{cor1}, 
is itself invariant.
  
Let $\phi$ be any bounded continuous function on $H^{-1}_0(\T)$.  
By invariance of $Q_{0, \b}^{(3)}$ under $\bar{\mathcal S}_t$, we have
\begin{equation*}
\int \phi \, dQ_{0, \b}^{(3)} = \int \phi\circ \bar{\mathcal S}_t \, dQ_{0, \b}^{(3)} .
\end{equation*}

\noi
Since $\bar{\mathcal S}_t$ is continuous on $H^{-1}_0(\T)$, 
we can take $\beta\to 0$ to obtain
\begin{equation*}
\int \phi \, dQ_0 = \int \phi\circ \bar{\mathcal S}_t \, dQ_0 
= \int \phi \, d \bar{\mathcal S}_t^*Q_0.
\end{equation*}

\noi
Taking $\phi(u)= \exp\big( i\jb{f, u} \big)$ for smooth mean zero functions $f$ on $\T$, 
we get
\begin{equation}
\int e^{ i\jb{f, u}}  d\bar{\mathcal S}_t^*Q_0 = e^{-\frac12\|f\|_{L^2}^2},
\end{equation}

\noi
which identifies $\bar{\mathcal S}_t^*Q_0 $ as mean zero white noise.  
This completes the straightforward proof of Theorem \ref{whiteinv}.  

The reason for calling the proof straightforward is that it is a fairly direct consequence of the intuitively
obvious fact  (\ref{conv1}).  It also has the advantage, partially exploited in the next subsection, that it does
not appear to rely on special properties of KdV.

\medskip

\begin{remark} \rm The same proof shows the invariance by KdV of mean zero white noise
$Q_{0,\sigma^2}$ with variance $\sigma^2$, defined
by
\begin{equation*}
\int e^{i \jb{f, u}} dQ_{0,\sigma^2}(u) = e^{-\frac{\sigma^2}2 \|f\|_{L^2}^2}.
\end{equation*}
\end{remark}

\subsection{Formal invariance of white noise for mKdV and cubic NLS on $\T$.}

\noi

The advantage of the straightforward proof of the invariance of white noise under the KdV flow presented in the previous
subsection is that it does not rely on special properties of KdV.  Hence, in principle, it provides a route towards invariance of
white noise for related equations. 

Unfortunately, Theorem \ref{thm1} is not enough to conclude the invariance of the white noise for mKdV or
cubic NLS (\eqref{NLS} with $p = 4$),
since their flows are  not expected to be well-defined below $H^{-\frac{1}{2}}$. Recall that 
mKdV and  cubic NLS are scaling-critical in $H^s$ with $s = -\frac{1}{2}$. 
{This means that 
the scaling invariance (on $\R$) 
$u(t,x)\mapsto \lambda^{-1} u(\lambda^{-2} t, \lambda^{-1} x)$ 
preserves the homogeneous $H^{-\frac{1}{2}}$-norm.}
It is usually expected that a nonlinear PDE is not well-posed below scaling-critical regularity, and the support of the white noise is 
below $H^{-\frac{1}{2}}$.   Nevertheless, if we lower our
standards, we are able to say something.  Let us define a measure $\mu$ to be {\it formally invariant} for a flow ${\mathcal S}_t$ if
there exist invariant measures $\mu_n$  for ${\mathcal S}_t$, converging weakly to $\mu$.

\begin{corollary} Mean zero white noise $Q_0$ is formally invariant for mKdV \eqref{mKdV}. 
\end{corollary}
\begin{corollary} Complex white noise ${\mathbf Q}$ is formally invariant for  cubic NLS (\eqref{NLS} with $p = 4$, either focussing or defocussing).
\end{corollary}

\begin{remark} \rm  Note that it is  {\it not} necessarily impossible to define the flows on the support of the white noise.
Indeed, one may be able to define the flow of mKdV or cubic NLS just on the support of the white noise.
See Bourgain \cite{B3}
for the case of the $L^2$-critical  defocusing cubic NLS on $\T^2$. 
The Gibbs measure on $\T^2$ is supported below $L^2 (\T^2)$.
Nonetheless, Bourgain constructed a well-defined flow on its support
and established the invariance of the Gibbs measure.  Also, given the formal invariance, it is very natural to expect that in these models, at least
${\mathcal S}^*_t$ has an extension to a class of measures including white noise.
\end{remark}

\begin{remark}  \rm The measures $Q^{(p)}_{0,\beta}$ are well defined for $2<p<6$, and all $\beta>0$.  Theorem \ref{thm1} extends 
readily to $2<p\le 4$.  $p=4$ is critical, in the sense that   $\beta\int_{\mathbb{T}} u^4 = \mathcal{O}(1)$ under
$Q^{(4)}_{0,\beta}$ as $\beta\to 0$, while for $2<p< 4$,  $\beta\int_{\mathbb{T}} u^p = o(1)$ under
$Q^{(p)}_{0,\beta}$.  For $p>4$, $\beta\int_{\mathbb{T}} u^p$ blows up.  Note that one should not conclude from this that 
Theorem \ref{thm1} cannot hold for $p>4$.  Indeed, it is quite plausible that it does.  However, the method of proof used here
does not extend beyond $p=4$.

\end{remark}

We conclude with some remarks on the concrete meaning of invariance vs formal invariance. 
Suppose that we want to start our dynamics, either KdV, mKdV, or cubic NLS, with $u_0$, distributed according to white noise.  
One way to proceed is to consider some
regularization $u^\beta_0 $, $\beta>0$, of the initial data $u_0$, 
and solve the equation in a more classical sense, to obtain smooth solutions
$u^\beta(t)={\mathcal S}_t u^\beta_0$ at a later time.  
Then, we ask if for small $\beta>0$, $u^\beta(t)$ is again approximately distributed
according to white noise.  
Invariance of white noise means that this procedure is true regardless of the type of regularization one uses.
Formal invariance means that there is at least one type of regularization which works:  In our case, the regularized $u^\beta_0 $ is distributed according to $Q_{0,\b}^{(4)}$.

\medskip

This paper is organized as follows.
In Section 2, we introduce the Wick-ordered monomials and prove a preliminary lemma.
In Section 3, we present the proof of Theorem \ref{thm1} for $p = 4$,
assuming the exponential expectation estimate \eqref{ExpEx1},
which we prove in Sections 4 and 5.
In Section 6, we briefly discuss the argument
for the complex-valued case, the defocusing case (Theorem \ref{thm0}),
and the $p = 3$ case.

\section{Wick ordering}

In this section, we perform a preliminary computation for the proof of Theorem \ref{thm1} for $p = 4$.
Recall that
\[d Q_{0, \b}^{(4)} = Z_\b^{-1} \, \ind_{\{\int_{\T} u^2 \leq K \b^{-\frac{1}{2}}\}}
e^{\b \int_{\T} u^4} d \mu_\b,\]

\noi
where $\mu_\b$ is as in \eqref{Gauss2}.
Under $\mu_\b$, $u$ is represented as a Fourier series \eqref{rep1},
where $g_n$ are independent standard complex Gaussians for $n>0$ and $g_{-n}=\overline g_{n}$. We will need various moments of $g_n$, the following identity can be proved e.g.~using the moment generating function of the complex Gaussian:
\begin{equation}\label{moments}
\mathbb{E} \left[g_n^k \,\overline{g_n^\ell}\right]=
 \dl_{k\ell} k!, \qquad k, \ell \in \Z_+,
\end{equation}

\noi
where $\dl_{k\ell} = 1$ if $k = \ell$ and $= 0$ otherwise.
In particular, $\mathbb{E} \left[g_{i_1} g_{i_2} \dots g_{i_k}\right]=0$ 
unless we can pair the indices $i_1,\dots, i_d$ in a way that the sum of the two indices is zero in each pair.

In order to study the behavior of $Q_{0, \b}^{(4)}$ as $\b \to 0$, 
we divide the space into several regions.
For this purpose, we introduce the Wick-ordered monomials 
$: u^2 \! :_\beta$ and $: u^4 \! :_\beta$ with parameter $\beta$:
\begin{align} \label{Wick2}
 & : u^2 \! :_\beta \, := u^2 -  a_\beta, \\
 \label{Wick4}
 & : u^4 \! :_\beta \, := u^4 - 6 a_\beta u^2 + 3 a_\beta^2,
\end{align}
\noindent
where
\[ a_\b = \mathbb{E}_{\mu_\b} \Big[\int_\T u^2 \Big] = \sum_{n \ne 0 } \frac{1}{1 + \betat n^2} .\]

\noi
For basics on Wick products and Gaussian Hilbert spaces, see e.g.~\cite{J}.
Note that $: u^k \! :_\b = H_k(u;a_\b)$, 
where $H(x, \s^2)$ is the Hermite polynomial in $x$ of degree $k$ with parameter $\s^2$.
We have 
\begin{equation}  \label{AA}
\b^\frac{1}{2} a_\beta \to  {\tfrac12} \quad \text{as }\ \beta \to 0,
\end{equation}

\noi
since 
$\beta^\frac{1}{2} \sum_{n \ne 0 } \frac{1}{1 + \betat n^2} \to  2 \int_{0}^\infty \frac{1}{1 +  4\pi^2 x^2} dx = \frac12  $
by Riemann sum approximation.
Also, by letting 
\[ b_\beta = \sum_{n \ne 0 } \frac{1}{(1 + \betat n^2)^2} 
\quad \text{ and } \quad c_\beta = \sum_{n \ne 0 } \frac{1}{(1 + \betat n^2)^4},\]

\noi
we have $\beta^\frac{1}{2}  b_\beta \to b_0$
and $\beta^\frac{1}{2}c_\beta \to c_0$ for some explicit constants $b_0, c_0 > 0$.

\begin{lemma} \label{LEM:expectation}
We have 
\begin{align}
\mathbb{E}_{\mu_\beta} \big[\int_\T :u^2 & \!  :_\beta \big] 
 = 0,    \quad
\mathbb{E}_{\mu_\beta} \Big[ \big( \int_\T :u^2 \!  :_\beta \big)^2 \Big] = 2 b_\beta,
\label{expectation1} \\
& \mathbb{E}_{\mu_\beta} \big[\int_\T :u^4 \!  :_\beta \big] = 0. \label{expectation2}
\end{align}

\noindent
Moreover, for sufficiently small $\beta>0$, we have\footnote{We 
use $A \lesssim B$ to denote an estimate of the form $A \leq CB$ for some $C >0$.
Similarly, we use $A \sim B$ to denote $A \lesssim B$ and $B \lesssim A$.} 
\begin{equation} \label{expectation3}
\mathbb{E}_{\mu_\beta} \Big[ \big( \int_{\T} :u^4 \!  :_\beta \big)^2 \Big]    \lesssim \beta^{-\frac{3}{2}} .
\end{equation}

\end{lemma}

\begin{proof}
For simplicity, we use $\mathbb{E}$ for $\mathbb{E}_{\mu_\beta}$.
By definition,  we have $\mathbb{E} [\int_\T u^2] = a_\beta$.
Also, we have 
\begin{align*}
\mathbb{E} \Big[ \big( \int_\T :u^2 \!  :_\beta \big)^2 \Big]
= 4 \mathbb{E} \bigg[\Big( \sum_{n\geq 1} \frac{|g_n|^2 -1}{ 1 + \betat n^2} \Big)^2\bigg]
= 4   \sum_{n\geq 1} \frac{\mathbb{E}[(|g_n|^2 -1)^2]}{ (1 + \betat n^2)^2} 
= 2 b_\beta.
\end{align*}

\noi
Using the representation of $u$ under $\mu_\beta$, we have 
\begin{equation} \label{L4}
 \int_\T u^4 = \sum_{\substack{n_{1234} = 0\\ \ n_j \ne 0}} 
\prod_{j = 1}^4 \frac{g_{n_j}}{\sqrt{ 1 +  \betat n_j^2 }},
\end{equation}

\noi
where $n_{1234} := n_1 + n_2 + n_3 + n_4$.
We say that we have a {\it pair} if we have $n_j = -n_k$, $j \ne k$ in the summation in \eqref{L4}.
Under the condition $n_{1234} = 0$, we have either two pairs or no pair.
Now, let $A_j = \{n_1 = -n_j\}$, $j = 2, 3, 4$,
Then, by symmetry, we can express the sum in \eqref{L4} as
\begin{align} \label{L4sum}
 \sum_{\substack{n_{1234} = 0\\ \ n_j \ne 0}} 
 = \sum_{\text{pair}} + \sum_{\text{no pair}} 
& = \sum_{j = 2}^4 \sum_{A_j} \ - \ \sum_{j<k} \sum_{A_j \cap A_k} \ + \ \sum_{\text{no pair}} \notag \\
& = 3 \sum_{\substack{n_1 = -n_2, \ n_3 = -n_4\\ n_j \ne 0 }} 
 -3 \sum_{\substack{n_1 = n_3 = -n_2 = -n_4\\ n_j \ne 0 }} 
 + \sum_{\text{no pair}} .
\end{align}

\noi
(Note that $A_2 \cap A_3 \cap A_4$ is empty.) 
From \eqref{Wick4}, we have
\begin{align} 
\int_\T :  u^4\!  :_\beta  \,  
& = 3 \Big\{ \sum_{n_1, n_3 \ne 0} \frac{|g_{n_1}|^2|g_{n_3}|^2}{(1 + \betat n_1^2)(1 + \betat n_3^2) }
- 2 a_\b \int_{\T} u^2 + a_\b^2 \Big\}  \notag \\
 &  \hphantom{XXX}
 - 3 \sum_{n\ne 0} \frac{|g_n|^4}{(1 + \betat n^2)^2}
 + \sum_{\text{no pair}} \prod_{j = 1}^4 \frac{g_{n_j}}{\sqrt{ 1 +  \betat n_j^2 }} \notag \\
& = 12 \bigg( \sum_{n_1 \geq 1} \frac{|g_{n_1}|^2 - 1}{1 + \betat n_1^2}\bigg)
\bigg(\sum_{n_3\geq 1} \frac{|g_{n_3}|^2 - 1}{1 + \betat n_3^2}\bigg) \notag \\
 &  \hphantom{XXX}
 - 6 \sum_{n \geq 1} \frac{|g_n|^4}{(1 + \betat n^2)^2}
 + \sum_{\text{no pair}} \prod_{j = 1}^4 \frac{g_{n_j}}{\sqrt{ 1 +  \betat n_j^2 }} \label{ZWick4}
 = : 12 \I_1 -6 \I_2 + \II.
\end{align}

\noi
Then, \eqref{expectation2} follows
from $\mathbb{E} [(|g_n|^2 -1)^2] = 1$ and $\mathbb{E}[|g_n|^4] ={2}$.
Using $\mathbb{E} [(|g_n|^2 -1)^4] = 9$, we have 
\begin{align*}
\mathbb{E} [\I_1^2] 
& =  \sum_{\substack{n_1, n_3 \geq 1 \\n_1 \ne n_3}} 
\frac{\mathbb{E}[(|g_{n_1}|^2 - 1)^2]}{(1 + \betat n_1^2)^2} 
\frac{\mathbb{E}[(|g_{n_3}|^2 - 1)^2]}{(1 + \betat n_3^2)^2}
+ \sum_{n \geq 1} \frac{\mathbb{E} [(|g_{n}|^2 - 1)^4]}{(1 + \betat n^2)^4} \\
& \leq \frac{b_\beta^2}{4} + \frac{9 c_\beta}{2} \lesssim \beta^{-1}
\end{align*}

\noi
for sufficiently small $\beta > 0$.
Similarly, we have
$\mathbb{E} [\I_2^2] \lesssim b_\b^2 +c_\b \lesssim  \beta^{-1}$.
Moreover, we have
\begin{equation} \label{ortho} 
\mathbb{E} [\I_1 \cdot \II] = \mathbb{E} [\I_2 \cdot \II] = 0.
\end{equation}
{by the comment after (\ref{moments}).}
Finally, we consider
\begin{align*}
\mathbb{E}[\II^2]
=  \mathbb{E} \bigg[\bigg( \sum_{\substack{n_{1234} = 0\\ n_j \ne 0\\ \text{no pair}}} 
\prod_{j = 1}^4 \frac{g_{n_j}}{\sqrt{ 1 +  \betat n_j^2 }} \bigg)
\bigg( \sum_{\substack{k_{1234} = 0\\ k_j \ne 0\\ \text{no pair}}} 
\prod_{j = 1}^4 \frac{g_{k_j}}{\sqrt{ 1 +  \beta k_j^2 }} \bigg)\bigg].
\end{align*}

\noi
{
Since the summation indices $\{ n_j\}$ and $\{ k_j\}$ contain no pair,
we see that the only nonzero contribution comes from $\{ n_1, n_2, n_3, n_4\} = - \{k_1, k_2, k_3, k_4\}$.} 
Thus, we have 
\begin{align*}
\mathbb{E}[\II^2]
=  24 \, \mathbb{E} \bigg[\bigg( \sum_{*} \prod_{j = 1}^4 \frac{|g_{n_j}|^2}{1 +  \betat n_j^2 } \bigg) \bigg]
\end{align*}

\noi
where $* = \{ n_{1234} = 0,  n_j \ne 0\text{, and no pair}\}$.
By separating the summation into (a) $n_j$ all distinct, (b) $n_1 = n_2 \ne n_3, n_4$ and $n_3 \ne n_4$, 
and (c) $n_1 = n_2 = n_3 \ne n_4$ (up to permutations of the indices),
we have 
\begin{align*}
\mathbb{E}[\II^2]
=  24 \bigg\{ \sum_{\substack{* \\ n_j,  \text{ all distinct}}} 
+ \, 6\cdot 2 \sum_{\substack{* \\ n_1 = n_2 \ne n_3, n_4 \\n_3 \ne n_4}} 
+ \, 4 \cdot 6 \sum_{\substack{* \\ n_1 = n_2 = n_3 \ne n_4}} \prod_{j = 1}^4 \frac{1}{1 +  \betat n_j^2 } \bigg\}
\end{align*}

\noi
since $\mathbb{E} [|g_n|^4] = 2$ and $\mathbb{E} [|g_n|^6] = 6$.
From the positivity of  the summands  and 
by Riemann sum approximation, we have 
\begin{align*}
 \mathbb{E}  [\II^2]  
 \lesssim &   \sum_{n_1, n_2, n_3 \ne 0} \prod_{j = 1}^3 \frac{1}{1 +  \betat n_j^2 } \frac{1}{1 + \betat (n_1 + n_2 +n_3)^2} \\
& + \sum_{n_1, n_3 \ne 0} \frac{1}{(1 +  \betat n_1^2)^2 } \frac{1}{1 +  \betat n_3^2 } \frac{1}{1 + \betat (2 n_1 +n_3)^2}  
+ \sum_{n_1, n_3 \ne 0} \frac{1}{(1 +  \betat n_1^2)^3 }  \frac{1}{1 + \betat (3 n_1)^2} \\
\sim  &  \  \betat^{-\frac{3}{2}} \int_{\R^3} \prod_{j = 1}^3 \frac{1}{1 +  x_j^2 } 
  \frac{1}{1 +  (x_1 + x_2 +x_3)^2} dx_1dx_2dx_3  \\
& + \betat^{-1} \int_{\R^2} \frac{1}{(1 +  x_1^2)^2 } \frac{1}{1 +   x_3^2 } 
\frac{1}{1 + (2 x_1 +x_3)^2} dx_1 dx_3  \\
& + \betat^{-\frac{1}{2}} \int_\R \frac{1}{(1 + x_1^2)^3 }  \frac{1}{1 + (3 x_1)^2} dx_1
 \lesssim  \beta^{-\frac{3}{2}}
\end{align*}

\noi
for sufficiently small $\beta >0$ .
Hence, we obtain  \eqref{expectation3}.
\end{proof}

\begin{remark} \label{REM:exp} \rm
The moral is that the main contribution of $\int_{\T} :  u^4\!  :_\beta$
comes from the ``no pair, all distinct'' part.
From \eqref{expectation2} and \eqref{AA}, we see that
$\mathbb{E}\big[\b \int_{\T} u^4\big] = 3 \b a_\b^2 = O(1).$
This shows that the decay of $\b$ and the growth of $\int_{\T} u^4$ is in perfect balance.
\end{remark}

\section{Proof of Theorem \ref{thm1}: $p = 4$ } 

 In order to prove Theorem \ref{thm1},
it suffices to show that, for any smooth mean 0 function $f$ on $\T$,  
\begin{equation} \label{mainconv}
C_\beta  \int  e^{i \int_{\T} f u + \beta \int_{\T} u^4 } \ind_{\{ \int_{\T} u^2 \leq K \beta^{-\frac{1}{2}}\}} d \mu_\beta 
\to  e^{-\frac{1}{2} \|f\|_{L^2}^2},
\quad \text{as} \quad \b \to 0
\end{equation}
for some $C_\beta >0$.  Indeed \eqref{mainconv} implies
\begin{equation}\label{wc}\int e^{i \int_{\T} f u} d Q_{0,\beta}^{(4)} =  
\frac{C_\beta \int  e^{i \int_{\T} f u + \beta \int_{\T} u^4 } 
\ind_{\{ \int_{\T} u^2 \leq K \beta^{-\frac{1}{2}}\}} d \mu_\beta }
{C_\beta \int  e^{ \beta \int_{\T} u^4 } 
\ind_{\{ \int_{\T} u^2 \leq K \beta^{-\frac{1}{2}}\}} d \mu_\beta } 
\to  \frac{ e^{-\frac{1}{2} \|f\|_{L^2}^2}}{e^{-\frac{1}{2} \|0\|_{L^2}^2}}
= e^{-\frac{1}{2} \|f\|_{L^2}^2}.
\end{equation}

\noi
This means that the joint distribution of the Fourier coefficients of $u$ 
under $Q_{0,\beta}^{(4)}$ converges weakly to the joint distribution of the coefficients from the white noise $Q_0$.  
The weak convergence of $Q_{0,\beta}^{(4)}$ to $Q_0$  in $H^s_0(\T)$, $s<-\frac12$, now follows from the following lemma, 
whose proof is presented at the end of this section. 
\begin{lemma}\label{LEM:tightness}
The sequence of measures  $Q_{0,\beta}^{(4)}$ is tight in $H^s_{0}(\T)$, $s<-\frac12$,
 as $\beta\to 0$.
\end{lemma}

\noi
{
It follows from from Lemma \ref{LEM:tightness} and Prohorov's theorem
that for any sequence $\{ \b_j\}$ of positive numbers tending to 0, 
the sequence $\big\{Q_{0,\beta_j}^{(4)}\big\}$ is sequentially compact.
Moreover, 
by the comment after \eqref{wc}, it converges weakly to $Q_0$.
The same comment guarantees the uniqueness of the limit point of 
$\big\{Q_{0,\beta}^{(4)}\big\}$ for $\beta \to 0$.
Hence, Theorem \ref{thm1} follows.}

\medskip

In view of Lemma \ref{LEM:expectation}, 
define $\mathcal{A}_{\beta, N}$ and $\mathcal{B}_{\beta, N}$ by
\begin{align} 
\mathcal{A}_{\beta, N}   = \big\{  \Big| \int_{\T} : u^4 \! :_\beta \Big| \leq N \beta^{-\frac{3}{4}}\big\},  
\text{ and }
\mathcal{B}_{\beta, N}  = \big\{  \Big| \int_{\T} :u^2 \!:_\beta \Big|  \leq N \beta^{-\frac{1}{4}} \big\} \label{setAB}
\end{align}

\noi
for large $N$ and small $\beta > 0$, and we consider separately the contributions from 
\[ \text{(i) }\mathcal{A}_{\beta, N} \cap \mathcal{B}_{\beta, N}, \quad
\text{ (ii) } \mathcal{A}_{\beta, N} \cap \mathcal{B}^c_{\beta, N}, \quad
\text{ and \quad (iii) } \mathcal{A}^c_{\beta, N}.\]

\noi
First, note that by Chebyshev's inequality with Lemma \ref{LEM:expectation} and \eqref{setAB}, 
we have an easy preliminary estimate
\begin{equation}\label{small}
 \mu_\beta\left( {\mathcal{A}_{\beta, N}^c \cup \mathcal{B}^c_{\beta, N}}\right) \lesssim N^{-2}.
\end{equation}

\noi
Our goal is to show that the main contribution for the weak convergence \eqref{mainconv} 
indeed comes from (i),
and that the contributions from (ii) and (iii) are small.

\medskip

\noi
$\bullet$ {\bf (i)} On  $\mathcal{A}_{\beta, N} \cap \mathcal{B}_{\beta, N}$: {
Since $\int_{\T} u^4=\int_{\T} :u^4\!:_\b+6 a_\b \int_{\T} u^2-3a_\b^2$ and 
 $\int_{\T} :u^4\!:_\b$ is ``small'' on $\mathcal{A}_{\beta, N} $, it is natural to introduce the 
 the Gaussian probability measure 
\begin{equation} \label{mubeta}
d \wt{\mu}_\beta = Z_\b^{-1}\exp \Big(6 \beta  a_\beta \int_\T  u^2\Big) \, d \mu_\beta
\end{equation}
for sufficiently small $\beta > 0$.}
First, we show that the normalization $Z_\b$ is indeed finite for (small) $\b > 0$.
\begin{lemma} \label{LEM:norm}
The normalization constant $Z_\b$ in (\ref{mubeta}) is bounded uniformly as $\beta \searrow 0$.
Moreover,
\[\lim_{\beta\to 0} \int e^{6 \beta  a_\beta \int_\T  u^2} \, d \mu_\beta
 = {e^{3/2}}\]
\end{lemma}

\begin{proof}
From \eqref{rep1}, we have, for small $\beta > 0$,  
\begin{align*}
\int e^{6 \beta  a_\beta \int_\T  u^2} \, d \mu_\b
& = \prod_{n\geq 1} \mathbb{E}\bigg[ \exp\Big( \frac{12 \b a_\b}{1+\betat n^2} |g_n|^2\Big)\bigg]
= \prod_{n\geq 1}\frac{1}{1-\frac{12 \b a_\b}{1+\betat n^2}}\\
& = \prod_{n\geq 1} \frac{1+ \frac{1}{\betat n^2}}{1+ \frac{1 - 12 \b a_\b}{\betat n^2}}
= \frac{\sinh (\pi \betat^{-\frac{1}{2}} )}{\pi \betat^{-\frac{1}{2}} }
\frac{\pi \sqrt{1 - 12 \b a_\b} \betat^{-\frac{1}{2}} }
{\sinh (\pi \sqrt{1 - 12 \b a_\b} \betat^{-\frac{1}{2}} )}.
\end{align*}

\noi
Here, we used $\mathbb{E}[ e^{aX^2} ] = (1-2a)^{-\frac{1}{2}}$, $a < \frac{1}{2}$,
for a {\it real}-valued standard Gaussian random variable $X$,
and the infinite product formula for $\sinh z$.
By \eqref{AA}, we have 
\[\lim_{\b\to 0} \int e^{6 \beta  a_\beta \int  u^2} \, d \mu_\beta
= \lim_{\b\to 0} \exp \big(\pi (\betat^{-\frac{1}{2}} - \sqrt{1 - 12 \beta a_\b} \betat^{-\frac{1}{2}})\big)
= e^{3/2}.\]
\end{proof}

Under $\wt{\mu}_\beta$, we have
\begin{equation} \label{representation2}
u(x) = \sum_{n \ne 0} \frac{g_n}{\sqrt{1 - 12 \beta  a_\beta + \betat n^2}} e^{2\pi i n x}.
\end{equation}

\noi
From \eqref{AA}, we have $ 12 \beta  a_\beta \sim \beta^\frac{1}{2} \to 0 $ as $\beta \to 0$,
{so this is well defined if $\beta$ is small enough.}
The following lemma, combined with the argument following \eqref{mainconv}, shows that  {the Fourier coefficients under $\wt{\mu}_\beta$ converge in distribution to those of the white noise}. 

\begin{lemma} \label{LEM:conv2}
There exists $C_\beta$, $\wt{C}_\beta> 0$ such that 
\begin{equation} \label{conv2}
\lim_{\beta \to 0}  C_\beta  \int e^{i \int_\T f u  + 6 \beta a_\beta \int_\T u^2 - 3 \beta a_\beta^2 } d \mu_\beta 
= \lim_{\beta \to 0}  \wt{C}_\beta  \int e^{i \int_\T f u  - 3 \beta a_\beta^2 } d \wt{\mu}_\beta 
=  e^{-\frac{1}{2} \|f\|_{L^2}^2}, 
\end{equation}

\noi
for any smooth mean 0 function $f$ on $\T$,
\end{lemma}

\begin{proof}By a direct computation, we have
\begin{align*}
\int e^{i \int_{\T} f u}  d \wt{\mu}_\beta 
& = \exp \Big\{ i  \sum_{n \ne 0} \frac{\ft{f}_n g_n }{\sqrt{1 - 12 \beta  a_\beta + \betat n^2}}  \Big\} \\
&  = \exp \Big\{ -\frac{1}{2}  \sum_{n \ne 0} \frac{|\ft{f}_n|^2 }{1 - 12 \beta  a_\beta + \betat n^2}  \Big\} 
\to e^{-\frac{1}{2} \|f \|_{L^2}^2}.
\end{align*}

\noi
Then, \eqref{conv2} follows from $e^{-3 \b a_\b^2} \to e^{{-3/4}}$ as $\beta \to 0$.
\end{proof}

Next, we show that $\beta \int_\T u^4$ is very close to $a_\beta \int_\T u^2$ in this case 
and that it does not affect the weak convergence in Lemma \ref{LEM:conv2}.
For conciseness of the presentation, let us define, for a function $F$ on $C(\T)$,
 \[I_f(F) = \int F(u)  e^{i \int_\T f u  + 6 \beta a_\beta \int_\T u^2 - 3 \beta a_\beta^2 } d \mu_\beta.\]

\begin{lemma} \label{LEM:bound3}
Let $K > {\frac12}$. Then, for $N>0$, we have
\begin{align} \label{bound3}
\limsup_{\beta \to 0} \bigg|\int_{\mathcal{A}_{\beta, N}\cap \mathcal{B}_{\beta, N}}  
\ind_{\{ \int_\T u^2 \leq K \beta^{-\frac{1}{2}}\}} & e^{i \int_\T f u +  \beta \int_\T u^4} d \mu_\beta 
- I_f(1)\bigg| 
\lesssim  N^{-1} .
\end{align}
\end{lemma}

\begin{proof}
On $\mathcal{A}_{\beta, N}$, we have 
$ \big|e^{\beta \int_\T : u^4  :_\beta } - 1\big|\lesssim \beta^\frac{1}{4} N $ for $\beta \leq  N^{-4}$.
Hence, we have 
\begin{align*}
 \bigg|\int_{\mathcal{A}_{\beta, N}\cap \mathcal{B}_{\beta, N}}  
 & \ind_{\{ \int_\T u^2 \leq K \beta^{-\frac{1}{2}}\}}  e^{i \int_\T f u +  \beta \int_\T u^4} d \mu_\beta 
  - I_f\big( \ind_{\mathcal{A}_{\beta, N}\cap \mathcal{B}_{\beta, N}}
  \ind_{ \{ \int_\T u^2 \leq K \beta^{-\frac{1}{2}}\}}\big) \bigg| \\
&  \leq    e^{ 6 \beta^\frac{1}{2} a_\beta K - 3\beta a_\beta^2}   
\int | e^{\beta \int_\T : u^4  :_\beta } - 1| d \mu_\beta
\lesssim  \beta^\frac{1}{4} N.
\end{align*}

\noi
since $6 \beta^\frac{1}{2} a_\beta K - 3\beta a_\beta^2 = O(1)$.
Moreover, on $\mathcal{B}_{\beta, N}$, 
given $\eps > 0$, there exists $\b_0>0$ such that  
\[ \int_\T u^2 = \int_\T : u^2 \! :_\beta + a_\b
\leq N \beta^{-\frac{1}{4}} + ({\tfrac12} + \tfrac{\eps}{2}) \beta^{-\frac{1}{2}}
\leq (\tfrac12 + \eps) \beta^{-\frac{1}{2}} \]

\noi
for $0< \b  < \b_0$.
Thus, we have $\mathcal{B}_{\beta, N} \subset \{ \int_{\T} u^2 \leq K \beta^{-\frac{1}{2}}\}$ 
for sufficiently small $\beta >0$ as long as $K > \frac12$.
Hence,  \eqref{bound3} follows once we show
\begin{equation}\label{bound31}
\limsup_{\beta \to 0}
|I_f( \ind_{\mathcal{A}_{\beta, N}\cap \mathcal{B}_{\beta, N}}) -I_f(1)| 
= \limsup_{\beta \to 0}
|I_f( \ind_{\mathcal{A}^c_{\beta, N}\cup \mathcal{B}^c_{\beta, N}})|
\lesssim N^{-1}.
\end{equation}

\noi
By Cauchy-Schwarz inequality along with \eqref{small}, we have
\begin{align} \label{ABC}
|I_f( \ind_{\mathcal{A}^c_{\beta, N}\cup \mathcal{B}^c_{\beta, N}})|
\leq \Big(\mu_\beta (\mathcal{A}^c_{\beta, N}\cup \mathcal{B}^c_{\beta, N})\Big)^\frac{1}{2}
\bigg(\int  e^{   6 \beta a_\beta \int_\T u^2  } d \mu_\beta\bigg)^\frac{1}{2}
\lesssim N^{-1}
\end{align}

\noi
since 
$\int e^{   6 \beta a_\beta \int_\T u^2  } d \mu_\beta = O(1)$ 
by Lemma \ref{LEM:norm}.
\end{proof}

\medskip
\noi
$\bullet$ {\bf (ii)} On $\mathcal{A}_{\beta, N} \cap \mathcal{B}^c_{\beta, N}$:
\quad
In this case, the Wick-ordered $L^4$-norm of $u$ is controlled.
Indeed, we have

\begin{lemma} \label{LEM:region2}
For $\beta \leq N^{-4}$, we have 
\begin{equation} \label{bound2}
\int_{\mathcal{A}_{\beta, N}\cap \mathcal{B}^c_{\beta, N}}  
\ind_{\{ \int_{\T} u^2 \leq K \beta^{-\frac{1}{2}}\}} e^{i \int_\T f u +  \beta \int_\T u^4} d \mu_\beta \lesssim N^{-2} .
\end{equation}

\end{lemma}

\begin{proof}
From \eqref{AA}, we have $\b^\frac{1}{2} a_\b = O(1)$.
Thus, on $ \mathcal{A}_{\beta, N}\cap\{ \int_\T u^2 \leq K \beta^{-\frac{1}{2}}\}$,
we have
\[ \beta \int_{\T} u^4 \leq \beta \Big|\int_{\T} :u^4 \!:_\beta\Big| + 6 \beta a_\beta \int_{\T} u^2 + 3 \beta a_\beta^2 
\lesssim 1 \]

\noi
for $\beta \leq N^{-4}$.
Then,  \eqref{bound2} follows from \eqref{small}.
\end{proof}

\medskip
\noi
$\bullet$ {\bf (iii)} On $\mathcal{A}^c_{\beta, N}$:
\quad In this case, we do not have any control on the the Wick-ordered $L^4$-norm of $u$.
Nonetheless, we have the following exponential expectation estimate.
\begin{proposition} \label{PROP:ExpEx1}
Let $ r> 0$. Then, we have
\begin{equation} \label{ExpEx2}
\mathbb{E}_{\mu_\beta} \big[\, \ind_{\{ \int_{\T} u^2 \leq K \b^{-\frac{1}{2}} \}} e^{r \b \int_{\T} u^4}\big]
= \int \ind_{\{ \int_{\T} u^2 \leq K \b^{-\frac{1}{2}} \}} e^{r \b \int_{\T} u^4} d \mu_\b
\leq C(r) <\infty,
\end{equation}

\noi
uniformly in small $\beta >0$.
\end{proposition}

\noi
For each {\it fixed} $\b > 0$, \eqref{ExpEx2} follows from \cite{LRS, B2}.
The difficulty lies in establishing the estimate uniformly in $\beta > 0$.
The proof requires both Fourier analytic and probabilistic approaches.
We present the proof of Proposition \ref{PROP:ExpEx1} in Sections 4 and 5.

\begin{lemma} \label{LEM:region3}
The following estimate holds uniformly in small $\beta > 0$.
\begin{equation} \label{region3}
\int_{\mathcal{A}^c_{\beta, N}}  \ind_{\{ \int_{\T} u^2 \leq K \beta^{-\frac{1}{2}}\}} 
e^{i \int_{\T} f u +  \beta \int_{\T} u^4} d \mu_\beta \lesssim N^{-1} .
\end{equation}
\end{lemma}

\begin{proof}
By the Cauchy-Schwarz inequality followed by \eqref{ExpEx2} and \eqref{small}, the left hand side of (\ref{region3}) is bounded by
\[
 \bigg(\int_{\mathcal{A}^c_{\beta, N}}  \ind_{\{ \int_{\T} u^2 \leq K \beta^{-\frac{1}{2}}\}}  d \mu_\beta\bigg)^{\frac{1}{2}} \bigg( \int \ind_{\{ \int_{\T} u^2 \leq K \b^{-\frac{1}{2}} \}} e^{2 \b \int_{\T} u^4} d \mu_\b\bigg)^{\frac{1}{2}} \lesssim N^{-1}.
\]
\end{proof}


Finally, \eqref{mainconv}
follows
from Lemmas \ref{LEM:conv2}, \ref{LEM:bound3}, \ref{LEM:region2}, \ref{LEM:region3}
by first taking $\b \to 0$ and then $N\to \infty$. Besides  proving  Proposition \ref{PROP:ExpEx1} (which is the content of the next two sections), the only part left is the proof Lemma \ref{LEM:tightness} which we present below.

\begin{proof}[Proof of Lemma \ref{LEM:tightness}]
For any measurable set $A$, we have
\begin{eqnarray} \notag
{Q_{0,\b}^{(4)}}(A)&=&\frac{\int \ind_A \, \ind_{\{ \int_{\T} u^2 \leq K \b^{-\frac{1}{2}} \}} e^{ \b \int_{\T} u^4} d \mu_\b}{\int  \ind_{\{ \int_{\T} u^2 \leq K \b^{-\frac{1}{2}} \}} e^{ \b \int_{\T} u^4} d \mu_\b} \le \frac{\big(\int_A d \mu_\b\big)^{\frac{1}{2}}\left(
\int  \ind_{\{ \int_{\T} u^2 \leq K \b^{-\frac{1}{2}} \}} e^{2 \b \int_{\T} u^4} d \mu_\b
\right)^{\frac{1}{2}}}{\int  \ind_{\{ \int_{\T} u^2 \leq K \b^{-\frac{1}{2}} \}} e^{ \b \int_{\T} u^4} d \mu_\b}\\
&\le& C \, \big\{\mu_\beta(A)\big\}^{\frac{1}{2}}.\label{eq:dom}
\end{eqnarray}

\noi
In the first line, we used the definition of $Q_{0,\b}^{(4)}$ and Cauchy-Schwarz inequality. 
The second line follows from Proposition \ref{PROP:ExpEx1} and from the fact that the denominator is bounded from below because of Chebyshev's inequality and Lemma \ref{LEM:expectation}:
\begin{equation}\label{lowerbound1}
\int  \ind_{\{ \int_{\T} u^2 \leq K \b^{-\frac{1}{2}} \}} e^{ \b \int_{\T} u^4} d \mu_\b
\geq K^{-1} \b^\frac{1}{2}\, \mathbb{E}_{\mu_\b} \Big[\int_\T u^2 \Big]
= K^{-1} \b^\frac{1}{2} a_\b \sim \tfrac{1}{2}K^{-1}>0.
\end{equation}

The upper bound \eqref{eq:dom} shows that 
it is enough to prove that the sequence $\mu_\b$ is tight in $H^s_0(\T)$ for $s=-\frac12-\eps$, $\eps>0$. 
Consider a probability space with the independent standard complex Gaussian random variables $g_n$ with $g_{-n}=\overline{g_n}$. Setting $u^{(\beta)}(x)=\sum_{n\neq 0} \frac{g_n}{\sqrt{1+\betat n^2}} e^{2\pi i n x}$ for $\beta\ge 0$,
 we have a joint realization of the measures $\mu_\b$ and $Q_0$. By the Borel-Cantelli lemma,
  we have $\sup_{n>0} \frac{|g_n|}{n^{\eps/2}}<\infty$ with probability one. 
This means that for the Fourier coefficients $\hat u_n^{(\beta)}$ of $u^{(\beta)}$, 
we have $|\hat u_n^{(\beta)}|\le C n^{\eps/2}$ a.s.~with a finite (but random) $C$. 
Since $\hat u_n^{(\beta)}\to \hat u_n^{(0)}=g_n$ a.s.~as $\b\to 0$ for all $n$, this implies that 
$u^{(\beta)}\to u^{(0)}$ a.s.~in $H^s_0(\T)$ for $s=-\frac12-\eps$. 
From Prohorov's theorem, 
we immediately have the tightness of the measures $\mu_\b$ and hence the statement of the lemma.
\end{proof}

%
%

\section{Bourgain's argument:  $\ld > \b^{-\frac{1}{2}-}$}

In this section and next, we present the proof of Proposition \ref{PROP:ExpEx1}.
It follows once we prove the following tail estimate.
\begin{lemma} \label{LEM:tail1}
There exists  $c, C>0$ and $\delta>0$ such that for all $\beta > 0$ and $\lambda \geq 1$, 
\begin{equation}\label{Tail1}
\mu_\beta \Big( \, \beta \|u\|_{L^4(\mathbb{T})}^4 > \ld, \int_{\T} u^2 \leq K \beta^{-\frac{1}{2}} \Big) 
\leq Ce^{-c \ld^{1 + \dl}}
\end{equation}
\end{lemma}

\noi
We will prove this lemma by considering two cases: 
$\ld>\beta^{-\frac12-}$ and $\ld<\beta^{-\frac12-}$.\footnote{We use 
$a+$ and $a-$ to denote $a+\eps$ and $a-\eps$, respectively, for arbitrarily small $\eps\ll 1$.}
For {\it fixed} $\beta > 0$, Bourgain \cite{B2} proved \eqref{Tail1}
via the dyadic pigeonhole principle with the large deviation estimate
(Lemma \ref{LEM:polar}.)
See Theorem \ref{THM:LRS} below.
In this section, we follow his approach to handle the case $\ld>\beta^{-\frac12-}$.
For this purpose, 
we need the following lemma on the tail probabilities of $\chi^2$ random variables.

\begin{lemma} \label{LEM:polar}
Let $g_1,g_2, \ldots$ be independent standard \emph{real}-valued Gaussian random variables.
Then for any $M\ge 1$, we have the following large deviation estimate:
\begin{equation} \label{polar}
P\bigg[\Big( \sum_{n=1}^M g_n^2\Big)^\frac{1}{2} \geq R \bigg]
\leq  e^{-\frac{1}{4}R^2}, \qquad R   \geq 3 M^\frac{1}{2} .
\end{equation}
\end{lemma}

\begin{proof}
By Markov's inequality, for $0\le t < 1/2$ we have
\[ P\bigg[\Big( \sum_{n=1}^M g_n^2\Big)^\frac{1}{2} \geq R \bigg]
\leq \frac{\mathbb{E}\big[\exp(t\sum_{n=1}^M g_n^2 )\big]}{\exp(t R^2)}=(1-2t)^{-\frac{M}{2}} e^{-t R^2}.\]

\noi
Choosing $t=\frac{1}{2}(1-\frac{M}{R^2})$, we get the upper bound 
\[\Big(\frac{R^2}{M}\Big)^\frac{M}{2}e^{-\frac{1}{2}R^2 + \frac{1}{2}M}
\leq e^{\frac{M}{2}\log(R^2/M) + \left(\frac{1}{18}-\frac{1}{2}\right) R^2}
\leq e^{-\frac{1}{4}R^2} \]

\noi
where in the last step we used that $\log x\le x/4$ for $x\ge 9$.
\end{proof}

Let us introduce some notations.
Given $M \in \mathbb{N}$, 
let $\proj_{>M}$ denote the Dirichlet projection onto the frequencies $\{|n| > M\}$.
i.e.  $\proj_{>M} u = \sum_{|n|>M} \ft{u}_n e^{2\pi i n x}.$
$\proj_{\leq M}$ is defined in a similar manner. 
Given $j \in \mathbb{N}$, let $M_j = 2^j M$.  
We use the notation $n\sim M_j$ to denote the set of integers $|n|\in (M_{j-1}, M_j]$,
and denote by $\proj_{M_j}$ the Dirichlet projection onto the dyadic block $(M_{j-1}, M_j]$
i.e.  
$\proj_{M_j} u = \sum_{n \sim M_j } \ft{u}_n e^{2\pi i n x}.$

\begin{lemma}\label{LEM:high} 
Let $p\geq 2$ and $\beta \leq 1$.  
Assume that 
$M\ge \max (\beta^{-\frac12-\delta}, \beta^{-\frac{p}2+1-\delta})$
for  some $\dl>0$.
Then there exists $c, C_1, C_2>0$ such that for $\ld \geq C_1$, 
\begin{equation} \label{highhigh}
\mu_\beta \left( \beta \| \proj_{ >M} u \|^p_{L^p(\mathbb T)} >\lambda\right) 
\leq C_2\exp\{ -c \lambda^{\frac2{p}} \beta^{ 1 - \frac2{p} } M^{ \frac2{p} + 1 }\}
\end{equation}
\end{lemma}

\begin{proof} 
Let $\sigma_j= C2^{-\epsilon j} $, $j=1,2,\ldots$ for some small $\epsilon>0$ where $C=C(\epsilon)$ is such that $\sum_{j=1}^\infty \sigma_j=1$.  
Then, we have
\begin{equation}\label{tth}
   \mu_\beta \Big(  \beta^{\frac1{p}}\| \proj_{ >M_0 } u \|_{L^p(\mathbb T)} >\lambda^{\frac1p} \Big) 
 \le   \sum_{j=0}^\infty \mu_\beta \left(\beta^{\frac1{p}}\| \proj_{ M_j } u \|_{L^p(\mathbb T)} >\sigma_j \lambda^{\frac1p}\right).
\end{equation}

\noi
There is a $c=c(p)<\infty$  such that for all $j=1,2,\ldots$, 
\begin{equation}
\| \proj_{M_j} u \|_{L^p(\mathbb T)} \le cM_j^{\frac12 - \frac1p}  \| \proj_{M_j} u \|_{L^2(\mathbb T)}. 
\end{equation}
This is the Sobolev inequality, though in this particular case it is a simple application of H\"older's inequality.

\noi
From \eqref{rep1}, we have 
$\| \proj_{M_j} u \|^2_{L^2(\mathbb T)} =\sum_{n\sim M_j} |\hat{u}_n|^2
=\sum_{n\sim M_j}  (1+\betat n^2)^{-1} |g_n|^2 $  
Hence, the right hand side of \eqref{tth} is bounded by 
\begin{equation}\label{67}
\sum_{j=0}^\infty P \bigg[(\sum_{n\sim M_j} g_n^2 )^{1/2} \ge R_j \bigg],\quad 
\text{where } R_j := \sigma_j \lambda^{\frac1p} \beta^{-\frac1{p}}M_j^{\frac1p - \frac12} 
(1+\beta M^2_j )^{1/2}.
\end{equation}

\noi
For $M\ge \max (\beta^{-\frac12-\delta}, \beta^{-\frac{p}2+1-\delta})$, we have
\begin{align*}
R_j \geq  C M^\eps  \ld^{\frac1p} \b^{\frac{1}{2}-\frac1{p}} M_j^{\frac1p + \frac12-\eps} 
\geq 3 M_j^\frac{1}{2}.
\end{align*}

\noi
By Lemma \ref{LEM:polar} to \eqref{67},
we conclude that \eqref{67} is bounded by
$\sum_{j=0}^\infty \exp\{-c  \sigma_j^2 \lambda^{\frac2p} \beta^{1-\frac2{p}}M_j^{\frac2p +1} \}$. 
This completes the proof.
\end{proof}

\noi
Before presenting the proof of Lemma \ref{LEM:tail1} for $\ld > \b^{- \frac{1}{2}-}$,
let us apply Lemma \ref{LEM:high} to prove the result in \cite{LRS, B2}.
Take $\beta = 1$, and let $\mu = \mu_1$.

\begin{theorem} 
[Lebowitz, Rose, and Speer \cite{LRS}, Bourgain \cite{B2}] \label{THM:LRS}  
Let $K<\infty$ and $r<\infty$.
For $2< p<6$,  and for $p=6$ with sufficiently small $K=K(r)>0$,
we have
\begin{equation}\label{LRS}
e^{ \int u^p } {\bf 1}_{\{ \int_{\T} u^2 \le K \} } \in L^r( d\mu).
\end{equation}
\end{theorem}

\begin{remark} \rm
The critical value $p=6$ is related to the $L^2$-criticality of the quintic NLS
and the quintic generalized KdV.
\end{remark}

\begin{proof}[Proof of Theorem \ref{THM:LRS}]
{
It is enough to prove that 
\[
\int_0^\infty  e^{\lambda} \, \mu\left(r \int u^p\ge \lambda, \int_{\T} u^2\le K\right) d \lambda < \infty.
\]
Let $M= c_0\lambda^{\frac2{p-2}} K^{-\frac{p}{p-2} }$
for some $c_0>0$.
By Sobolev inequality, 
\[ \|\proj_{\le M} u \|_{L^p(\mathbb T)}
\le cM^{\frac12 - \frac1{p}} \|\proj_{\le M} u \|_{L^2(\mathbb T)}.\]
Hence,  we have $r\|\proj_{\le M} u \|^p_{L^p(\mathbb T)}\le  \ld/2$ on $\int_{\T} u^2 \leq K$.
For sufficiently large $\lambda >0$, 
the condition of  Lemma \ref{LEM:high} holds, so we have
\begin{equation}\label{LRSbnd}
\mu \Big(  r \| \proj_{ >M } u \|^p_{L^p(\mathbb T)} >\ld\Big) 
\le C\exp\{ -c r^{-\frac{p}{2}}\lambda^{\frac2{p}}  M^{ \frac2{p} + 1 }\}
= C \exp\{ -c' \lambda^{1+ \frac{6-p}{p-2}} r^{-\frac{p}{2}} K^{-\frac{p+2}{p-2}} \}.
\end{equation}
and the statement follows. Note that when $p = 6$, we need to take $K = K(r)$ sufficiently small
such that $r^{-3} K^{-2}$ is large and the coefficient of $\lambda$ is less than $-1$ in (\ref{LRSbnd}).}
\end{proof}

Now, we present the proof of Lemma \ref{LEM:tail1} for $ \ld > \beta^{-\frac{1}{2}-}$.
As we see, one obtains much less
in estimating the tail uniformly in $\beta>0$ even when $p = 4$.   
Indeed, Bourgain's argument is not enough to conclude the argument even for $p = 3$.

\begin{proof} [Proof of Lemma \ref{LEM:tail1} for $\ld > \beta^{-\frac{1}{2}-}$] 
The proof is similar to that of Theorem \ref{THM:LRS}.
First, choose $M = c_0K^{-2}\lambda >\beta^{-\frac{1}{2}-} $.  
By Sobolev inequality, 
\[\beta \|\proj_{\le M} u \|_{L^4(\mathbb T)}^4
\le c \beta M \|\proj_{\le M} u \|_{L^2(\mathbb T)}^4.\]  

\noi
Hence, on $\|\proj_{\le M} u \|_{L^2(\mathbb T)} \le K^{\frac{1}{2}}\beta^{-\frac{1}{4}}$, 
we have, for sufficiently small $c_0$,
 \[\beta\|\proj_{\le M} u \|^4_{L^4(\mathbb T)}\le c \, c_0^4 \ld \leq \lambda/2.\]

\noi
As before, we can apply Lemma \ref{LEM:polar} to handle the high frequencies
as long as $R_j \geq 3 M_j^{\frac{1}{2}}$ in \eqref{67}.
Unlike the proof of Lemma \ref{LEM:high}, when checking this,
we use the non-smallness of $M_j \geq\ld > \beta^{-\frac{1}{2}-}$.
In this case, we have
\begin{align*}
R_j & = \sigma_j \lambda^{\frac{1}{4}} \beta^{-\frac{1}{4}}M_j^{-\frac{1}{4}} 
(1+\beta M^2_j )^{1/2} \geq \beta^{\frac{1}{8} - } M_j^{\frac{3}{4}-\eps}
\geq M_j^{\frac{1}{2}+}.
\end{align*}

\noi
By proceeding as in the proof of Lemma \ref{LEM:high}, 
we obtain \eqref{highhigh}.
Then, \eqref{Tail1} follows
once we note that $M_j \gtrsim \ld  > \beta^{-\frac{1}{2}-}$.
\end{proof}

\section{Hypercontractivity estimate: $\ld < \b^{-\frac{1}{2}-}$}
{First, note that we have
$\beta \int_{\T} u^4 = \beta \int \! : \! u^4 \!:_\beta + O(1)$ on $\{ \int_{\T} u^2 \leq K \beta^{-\frac{1}{2}}\}$ and thus 
it is enough to prove \eqref{Tail1} with $\beta \int \! : \! u^4 \!:_\beta$ instead of $\beta \int_{\T} u^4$. We will use the identity \eqref{ZWick4} and we further separate
 the summation for $\II$ into (a) $n_j$ all distinct, (b) $n_1 = n_2 \ne n_3, n_4$ and $n_3 \ne n_4$, 
and (c) $n_1 = n_2 = n_3 \ne n_4$ (up to permutations of the indices)
and write $\II = \II_a + \II_b + \II_c$.  Recall also the definitions of $\I_1$ and $\I_2$ from  \eqref{ZWick4}.
We will show that the main contribution of $\beta \int \! : \! u^4 \!:_\beta$ comes from ``no pair, all distinct", i.e.~$\II_a$.}
\begin{lemma} \label{lemma62} 
On $\{ \int_{\T} u^2 \leq K \beta^{-\frac{1}{2}}\}$,  
there is a $C<\infty$ such that 
$\beta |\I_1|, \beta |\I_2|, \beta|\II_b|, \beta|\II_c|\le C$
uniformly in $\beta > 0$.
\end{lemma}

\begin{proof}
In view of \eqref{rep1}, we have 
\begin{align*}
\beta|\I_1| = \beta \bigg( \sum_{n\geq1} \frac{|g_n|^2-1 }{1+ \betat n^2} \bigg)^2 
\leq 2\beta \bigg( \sum_{n \geq 1} \frac{|g_n|^2}{1+\betat n^2} \bigg)^2
+ 2\beta \bigg( \sum_{n \geq 1} \frac{1}{1+\betat n^2}\bigg)^2 
\lesssim 1 
\end{align*}

\noi
on $\{ \int_{\T} u^2 \leq K \beta^{-\frac{1}{2}}\}$.
By H\"older inequality and $l^2 \subset l^4$, 
the contribution for $\II$ from the case (c) is at most
\begin{align*}
\beta|\II_c| & \sim \beta \bigg|\sum_{n_1 \ne 0} 
\frac{g_{n_1}^3}{(1+ \betat n_1^2)^\frac{3}{2}} \frac{g_{-3n_1} }{\sqrt{1 + \betat (-3n_1)^2}}\bigg|
\leq \beta \sum_{n \ne 0} \frac{|g_n|^4}{(1+\betat n^2)^2} \\
& \leq \beta \bigg(\sum_{n \ne 0} \frac{|g_n|^2}{1+\betat n^2}\bigg)^2 
= \beta \bigg( \int_{\T} u^2 \bigg)^2 \lesssim 1.
\end{align*}

\noi
Similarly, we have $\beta|\I_2| \lesssim 1$. 
Then, 
the contribution for $\II$ from the case (b) is at most
\begin{align*}
\beta|\II_b| & \sim \beta \bigg|\sum_{\substack{\text{no pair}\\ n_1, n_3 \ne 0}} 
\frac{g_{n_1}^2}{1+ \betat n_1^2} \frac{g_{n_3} }{\sqrt{1 + \betat n_3^2}}
\frac{g_{-2n_1-n_3} }{\sqrt{1 + \betat (-2n_1 -n_3)^2}}\bigg| \\
&\leq \beta \sum_{n_1 \ne 0} \frac{|g_{n_1}|^2}{1+\betat n_1^2} 
\sup_{n_1 \ne 0} \bigg|\sum_{n_3\ne 0} \frac{g_{n_3} }{\sqrt{1 + \betat n_3^2}}
\frac{g_{-2n_1-n_3} }{\sqrt{1 + \betat (-2n_1 -n_3)^2}}\bigg| \\
& \leq \beta \bigg(\sum_{n \ne 0} \frac{|g_n|^2}{1+\betat n^2}\bigg)^2 
= \beta \bigg( \int_{\T} u^2 \bigg)^2 \lesssim 1,
\end{align*}

\noi
where we used $ab \leq a^2/2 + b^2/2$ in the last line.
\end{proof}

In estimating the contribution from $\II_a=$``no pair, all distinct'',
we will use the hypercontractivity of the Ornstein-Uhlenbeck process.
Let $L$ denote the generator of the Ornstein-Uhlenbeck process on $H := L^2(\R^d, e^{-|x|^2/2}dx)$
given by $L = \Dl - x \cdot \nabla$.
Then, let $S(t) = \exp( tL)$
be the semigroup associated with $\dt u = L u$.
Then, the hypercontractivity of the Ornstein-Uhlenbeck semigroup \cite[Sec.3]{T}
says the following:
\begin{lemma} \label{LEM:HYP1}
Let $q \geq 2$. For $f \in H$ and $t\geq \frac{1}{2} \log(q-1)$,
we have
\[ \| S(t) f \|_{L^q(\R^d, \exp(-|x|^2/2)dx)}
\leq \|  f \|_{L^2(\R^d, \exp(-|x|^2/2)dx)}\]
\end{lemma}

The eigenfunctions of $L$ are given by $\prod_{j = 1}^d h_{k_j}(x_j)$, 
where $h_k$ is the Hermite polynomial of degree $k$,
and the corresponding eigenvalue is given by
 $\ld = - (k_1+ \cdots + k_d)$.
The first few Hermite polynomials are
\begin{equation*} 
h_0(x) = 1,  \ h_1(x) = x, \ h_2(x) = x^2 - 1, \ \ldots
\end{equation*}

\noi
Let \[H(x) = \sum_\G c(n_1, \ldots, n_4) x_{n_1} \cdots x_{n_4},\] 

\noi
where 
$\G = \{ (n_1, \cdots, n_4) \in \{1, \cdots, d\}^4, \text{ all distinct} \}$.
Note that $H(x)$ is an eigenfunction of $L$ with the eigenvalue $-4$.  
The following {\it dimension-independent} estimate is a simple consequence of Lemma \ref{LEM:HYP1}:

\begin{corollary}  
For all $d=1,2,3,\ldots$, we have
\begin{equation} \label{LEM:HYP2}
 \| H(x) \|_{L^{q}(\R^d, \exp(-|x|^2/2)dx)}
\leq q^2 \|  H(x)  \|_{L^2(\R^d, \exp(-|x|^2/2)dx)}.
\end{equation}
\end{corollary}

\begin{proof} [Proof of Lemma \ref{LEM:tail1} for $\ld < \beta^{-\frac{1}{2}-}$] 
By Lemma \ref{lemma62} and the argument just preceding it, all 
it suffices to prove
\begin{equation}\label{Tail1.5}
\mu_\beta \Big( |\II_a| \ge \ld, \int_{\T} u^2 \leq K \beta^{-\frac{1}{2}} \Big) 
\leq C e^{-c \ld^{1 + \dl}}
\end{equation}

\noi
for $\ld \leq \beta^{-\frac{1}{2}-}$.
First, we show
\begin{equation}\label{Tail2}
\mu_\beta \Big( |F_{\beta,M}| \ge \ld, \int_{\T} u^2 \leq K \beta^{-\frac{1}{2}} \Big) 
\leq C e^{-c \ld^{1 + \dl}}
\end{equation}

\noi
for $\ld \leq \beta^{-\frac{1}{2}-}$,  where
\begin{equation} \label{QQ}
F_{\beta,M} = \beta \sum_{**} \prod_{j = 1}^4 \frac{g_{n_j}}{\sqrt{ 1 +  \betat n_j^2 }} 
\end{equation}
with
$** = \{ \, n_{1234} := n_1 + \cdots + n_4 = 0,  n_j \ne 0,  
\text{ no pair, all distinct, } |n_j|\leq M\}$, with a constant $c$ independent of $M$.
Then, we will indicate how \eqref{Tail1.5} follows from \eqref{Tail2}.

By expanding the complex-valued Gaussians $g_n$ into their real and imaginary parts, 
we can apply \eqref{LEM:HYP2} to $Q_{\beta,M} $ in \eqref{QQ}.
From (the proof of) Lemma \ref{LEM:expectation}, we have
$\| F_{\beta,M} \|_{L^2(d\mu_\beta)} \leq C \beta^{\frac{1}{4}}$.
By \eqref{LEM:HYP2}, 
we have \begin{equation} \label{QQQ}
\| F_{\beta,M} \|_{L^q(d\mu_\beta)} 
\leq C q^2 \beta^\frac{1}{4}
\end{equation}

\noi
for all $ q\geq 2$.
Note that we need that $u$ has a finite Fourier support, but 
the actual upperbound on the support is not important.
Then, we have
\begin{equation}\label{QQQQ}
 \int \exp( c \beta^{-\frac{1}{8}} | F_{\beta,M} |^\frac{1}{2} ) d\mu_\beta \leq C
 \end{equation}

\noi 
from Lemma 4.5 in \cite{T}. This can be proved 
by expanding the exponential in the Taylor series and applying \eqref{QQQ} and H\"older's inequality. 
Equation \eqref{QQQQ} in turn implies
$ \mu_\beta ( | F_{\beta,M} | > \ld ) \leq C \exp( -c' \beta^{-\frac{1}{8}} \ld^\frac{1}{2} )$
by Markov's inequality, 
i.e.
we proved \eqref{Tail2} for $\ld \leq \beta^{-\frac{1}{4}+}$.

Now, we consider the remaining case: $ \beta^{-\frac{1}{4}+} \leq \ld \leq \beta^{-\frac{1}{2}-}$.
Then, using $\ld \geq \beta^{-\frac{1}{4}+\eps}$,
\begin{align*}
\mu_\beta \big( | F_{\beta,M} | \geq \ld ) & \leq \frac{\| F_{\beta,M} \|^q_{L^q(d \mu_\beta)}}{\ld^q}
\leq C q^{2q} \beta^{\frac{q}{2} -\eps q}
\leq e^{2q \ln q} e^{-\frac{q}{3} \ln \beta^{-1}}
= e^{ - \frac{q}{3} \ln \beta^{-1} + 2q \ln q}
\intertext{By choosing $q \sim \beta^{-\frac{3}{4}} \ll \beta^{-1}$ and using $\ld \leq  \beta^{-\frac{1}{2}-\eps}$, }
& \leq e^{ - c \beta^{-\frac{3}{4}} \ln \beta^{-1}}
\leq e^{-c\ld^{\frac{3}{2}-}}.
\end{align*}

\noi
This proves \eqref{Tail2}.

Now, we need to show how \eqref{Tail1.5} follows from \eqref{Tail2}.
Clearly, $F_{\b,M}\to \II_a$ in $L^2(d \mu_\b)$ as $M\to \infty$. Thus, we can find a subsequence $M_k\to\infty$ for which $F_{\b,M_k}\to \II_a$ almost surely with respect to $\mu_\b$. 
By the dominated convergence theorem for the indicator random variables $\ind\left( |F_{\beta,M_k}| \ge \ld, \int_{\T} u^2 \leq K \beta^{-\frac{1}{2}}  \right)$, we have, for fixed $\beta >0$ and $\ld \geq 1$,  
\[
\mu_\beta \Big( |\II_a| \ge \ld, \int_{\T} u^2 \leq K \beta^{-\frac{1}{2}} \Big) 
=\lim_{k\to\infty} \mu_\beta \Big( |F_{\b, M_k}| \ge \ld, \int_{\T} u^2 \leq K \beta^{-\frac{1}{2}} \Big) \le C e^{-c \ld^{1 + \dl}},
\]

\noi
where $C$ and $c$ are independent of $\beta$ and $\ld$.
This completes the proof of the tail estimate \eqref{Tail1}.
\end{proof}

\section{Remarks}

We proved Theorem \ref{thm1} for $p = 4$.
In this section, we briefly discuss the minor changes needed to handle
the complex-valued case, the focusing case (Theorem \ref{thm0}),
and the $p = 3$ case.

\medskip

\noi
$\bullet$ {\bf Complex-valued case:}
As mentioned in Remark \ref{REM:complex},
the same result holds for the complex-valued case as well.
In this case, one needs to use the following definitions of Wick-ordered monomials,
\begin{align*}
 & : |u|^2 \! :_\beta \, = |u|^2 -  a_\beta, \\
  & : |u|^4 \! :_\beta \, = |u|^4 - 4 a_\beta |u|^2 + 2 a_\beta^2,
\end{align*}

\noindent
where $ a_\b = \mathbb{E}_{\mu_\b} \big[\int_\T |u|^2 \big].$
The proof is basically the same (note that we did not really need the mean-zero condition),
and one needs to prove Proposition \ref{PROP:ExpEx1} in the complex-valued case.
This follows easily once we note $|u|^4 \lesssim (\text{Re} \, u)^4 +  (\text{Im} \, u)^4$.

\medskip

\noi
$\bullet$ {\bf Defocusing case:}
Now, let us briefly discuss the proof of Theorem \ref{thm0}.
First, write 
\begin{align*}
\int e^{i \int fu} & d \wt{Q}_{0, \b}^{(4)} 
= Z_\b^{-1} \int e^{i \int fu - \beta \int_{\T} u^4} d \mu_\b\\
& = Z_\b^{-1} \int e^{i \int fu - \beta \int_{\T} u^4} \ind_{\{\int_{\T} u^2 \leq K \b^{-\frac{1}{2}}\}} d \mu_\b
+ Z_\b^{-1} \int e^{i \int fu - \beta \int_{\T} u^4} \ind_{\{\int_{\T} u^2 > K \b^{-\frac{1}{2}}\}} d \mu_\b.
\end{align*}

\noi
By repeating the argument in Section 3, 
the first term yields the desired result.
Note that we have Proposition \ref{PROP:ExpEx1} for free thanks to the negative sign.
As for the second term, \eqref{small} states that the contribution on 
$\mathcal{A}_{\beta, N}^c \cup \mathcal{B}^c_{\beta, N}$
goes to 0 as $N\to \infty$.
The contribution on $\mathcal{A}_{\beta, N} \cap \mathcal{B}_{\beta, N}$
also goes to 0 since
$\mathcal{A}_{\beta, N} \cap \mathcal{B}_{\beta, N}
\subset
\mathcal{B}_{\beta, N} \subset \{ \int_{\T} u^2 \leq K \beta^{-\frac{1}{2}}\}$ 
for sufficiently small $\beta >0$ for $K>\frac12$.

Note that Lemma \ref{LEM:tightness} follows
in a similar manner as before, 
once we show that 
the denominator in \eqref{eq:dom} 
is bounded from below.
%
%
%
%
%
%
By Jensen's inequality we have
\begin{equation}\label{eq:lower}
\int_{\mathcal{A}} e^{ -\b \int_{\T} u^4} d \mu_\b
\geq\mu_\b(\mathcal{A}) \exp\left\{ - \frac{1}{\mu_\b(\mathcal{A})}   \mathbb{E}_{\mu_\b}\left[\ind_{\mathcal{A}} \b\int_{\T} u^4 \right]    \right\}
\end{equation}
where $\mathcal{A}=\{ \int_{\T} u^2  \leq K \b^{-\frac{1}{2}} \}$. The right hand side is clearly bounded from below as $\b\to0$ since $\b \int_\T u^4\to C$ by Lemma \ref{LEM:expectation} and $\mu_\b(\mathcal{A}) $ is bounded from below by Chebyshev (c.f. (\ref{lowerbound1})).

\medskip

\noi
$\bullet$ {\bf $p = 3$ case:}
The proof of Theorem \ref{thm1} for $p = 3$ 
is similar to the $p = 4$ case.  Once we have Lemma \ref{LEM:tail1}, everything follows for $p<4$.
However, in this case, we do not need to use the Wick-ordered $\int_\T u^3$, and a simpler proof is 
available because the hypercontractivity estimates
can be replaced by a direct application of the Sobolev inequality, but it is still  a nontrivial extension of the Bourgain
method.  We sketch it now.

By direct computation, we have
\[\mathbb{E}_{\mu_\b}\bigg[\int_{\T} u^3\bigg]  =  0,  \quad \text{ and }\quad
\mathbb{E}_{\mu_{\b}}\bigg[\Big(\int_{\T} u^3\Big)^2\bigg] \lesssim \b^{-1}.\]

\noi
Similarly to the $p=4$ case  we define $\mathcal{C}_{\beta, N}$ by
\begin{align} 
\mathcal{C}_{\beta, N}  = \big\{  \Big| \int_{\T} u^3 \Big|  \leq N \beta^{-\frac{1}{2}} \big\} \label{setAC},
\end{align}

\noi
and
separately estimate the contributions from 
\[ \text{(i) }\mathcal{A}_{\beta, N} \cap \mathcal{C}_{\beta, N}, \quad
\text{ (ii) } \mathcal{A}_{\beta, N} \cap \mathcal{C}^c_{\beta, N}, \quad
\text{ and \quad (iii) } \mathcal{A}^c_{\beta, N}.\]

\noi
The main contribution comes from $\mathcal{A}_{\beta, N} \cap \mathcal{C}_{\beta, N}$.
Unlike the $p = 4$ case, there is no need to introduce $\wt{\mu}_\b$ defined in \eqref{mubeta},
and 
we can simply use the convergence of $\mu_\b$:
$\lim_{\b\to 0} \int e^{ i \int_{\T} f u}d\mu_\b = e^{-\frac12\|f\|_2^2}$
for any mean zero smooth function $f$ on $\T$.

The contributions from  $\mathcal{A}_{\beta, N} \cap \mathcal{C}^c_{\beta, N}$
and $\mathcal{A}^c_{\beta, N}$
can be shown to be small by Chebyshev's inequality, 
once we prove the following exponential expectation bound.

\begin{proposition} \label{PROP:Ex3}
Let $ r> 0$. Then, we have
\begin{equation} \label{Ex3}
\mathbb{E}_{\mu_\beta} \big[\, \ind_{\{ \int_{\T} u^2 \leq K \b^{-\frac{1}{2}} \}} e^{r \b \int_{\T} u^3}\big]
= \int \ind_{\{ \int_{\T} u^2 \leq K \b^{-\frac{1}{2}} \}} e^{r \b \int_{\T} u^3} d \mu_\b
\leq C(r) <\infty,
\end{equation}

\noi
uniformly in small $\beta >0$.
\end{proposition}

\noi
Proposition \ref{PROP:Ex3} is a corollary 
of Proposition \ref{PROP:ExpEx1}.
However, there is an easier direct proof in this case:

\begin{proof}
By Sobolev inequality followed by H\"older's inequality, we have 
\begin{align*}
\int_\T u^3  & \leq c \bigg(\sum_{n \ne 0} n^\frac{1}{3} |\hat u_n|^2 \bigg)^{\frac{3}{2}}
 \leq  c\bigg(\sum_{n\ne 0} n^\frac{1}{2} |\hat u_n|^2\bigg) \bigg(\sum_{n\ne0} |\hat u_n|^2\bigg)^{1/2} \leq cK^\frac{1}{2} \b^{-\frac{1}{4}} \sum_{n\ne0} n^\frac{1}{2} |\hat u_n|^2
\end{align*}

\noi
on $\int_\T u^2 \leq K \b^{-\frac{1}{2}}$.
Moreover, we have
\[ \beta \int_\T \big[\, \proj_{\leq c_0 \b^{-\frac{1}{2} }}u\big]^3
\lesssim
\b^\frac{3}{4} \sum_{1 \leq |n| \leq c_0 \b^{-\frac{1}{2} }}  n^\frac{1}{2} |\hat u_n|^2
\leq C\]

\noi
on $\int_{\T} u^2 \leq K \b^{-\frac{1}{2}}$.
Hence, from \eqref{rep1}, we have
\begin{align} \label{ex31}
\int \ind_{\{ \int_\T u^2\le K\b^{-\frac{1}{2}} \} } & e^{ r\beta\int_{\T} u^3 } d\mu_\b
 \leq \int \ind_{ \{\int_{\T} u^2\le K\b^{-\frac{1}{2}}\}}
\exp\Big\{ C + c\b^\frac{3}{4} \sum_{|n| > \, c_0\b^{-\frac{1}{2}}}  n^\frac{1}{2} |\hat u_n|^2 \Big\} d\mu_\b \notag\\
& \lesssim \int \prod_{n > \, c_0 \b^{-\frac{1}{2}}}
\exp\bigg\{  \frac{2c\b^\frac{3}{4} n^\frac{1}{2}}{1+\betat n^2} |g_n|^2 \bigg\} d\mu_\b
 =  \prod_{n > \, c_0\b^{-\frac{1}{2}}} \frac{1}{1-\frac{c\b^\frac{3}{4} n^\frac{1}{2}}{1+\betat n^2}},
\end{align}

\noi
where in the last equality we used $\mathbb{E}[ e^{aX^2} ] = (1-2a)^{-\frac{1}{2}}$, $a < \frac{1}{2}$
for a  real valued standard Gaussian random variable $X$,
since $(c\b^\frac{3}{4} n^\frac{1}{2})(1+\betat n^2)^{-1} <\frac{1}{2}$
on $n > \, c_0\b^{-\frac{1}{2}}$ for sufficiently large $c_0>0$.

It is not hard to check that  $0<x<1/2$ implies $(1-x)^{-1}<e^{x+ x^2}$.  
\begin{align*}
\eqref{ex31} & \leq
 \prod_{n > \, c_0\b^{-\frac{1}{2}}} \exp\bigg\{
 \frac{c\b^\frac{3}{4} n^\frac{1}{2}}{1+\betat n^2} +
 \frac{c^2 \b^\frac{3}{2} n}{(1+\betat n^2)^2}\bigg\} =  \exp\bigg\{
 \sum_{n > \, c_0\b^{-\frac{1}{2}}}\frac{c\b^\frac{3}{4} n^\frac{1}{2}}{1+\betat n^2} +
 \frac{c^2 \b^\frac{3}{2} n}{(1+\betat n^2)^2}\bigg\}.
\end{align*}

\noi
Hence, by Riemann sum approximation, we have for sufficiently small $\beta > 0$, 
\begin{align*}
& \sum_{n> c_0 \,  \beta^{-\frac{1}{2}}}  \frac{\beta^\frac{3}{4} n^\frac{1}{2}}{1 +\betat n^2} + \frac{\beta^{\frac{1}{2}} n}{(1 +\betat n^2)^2}
\lesssim \int_{c_0}^\infty \frac{\sqrt{x}}{1+x^2}dx+\b^{1/2}\int_{c_0}^\infty \frac{x}{(1+x^2)^2}dx <\infty.
\end{align*}
%
\noi
This shows that \eqref{ex31} is finite.
\end{proof}

%
Lastly, note that Lemma \ref{LEM:tightness} follows as before, 
once we show that 
the denominator in \eqref{eq:dom} (with $p = 3$)
is bounded from below.
Proceeding the same way as in \eqref{eq:lower} this is immediate since $\mathbb{E}_{\mu_\b}\left[ \ind_{{\int_{\T} u^2  \leq K \b^{-\frac{1}{2}} }} \b \int_{\T} u^3 \right]  =0$ by the $u \to -u$ symmetry of $\mu_\b$.

\end{document}